\documentclass{amsart}
\usepackage{amsfonts,amssymb,amsmath,amsthm}
\usepackage{url}
\usepackage{enumerate}
\usepackage{graphicx}

\newtheorem{thrm}{Theorem}[section]

\newtheorem{prop}[thrm]{Proposition}
\newtheorem{cor}[thrm]{Corollary}
\theoremstyle{definition}
\newtheorem{definition}[thrm]{Definition}

\numberwithin{equation}{section}

\author{Don Colladay}
\address{
Division of Natural Science\\
New College of Florida\\
Sarasota, FL}
\email{colladay@ncf.edu}

\author{Leon Kaganovskiy}
\address{
Mathematics Department\\
Touro College\\
Brooklyn, NY}
\email{leonkag@gmail.com}

\author{Patrick McDonald}
\address{
Division of Natural Science\\
New College of Florida\\
Sarasota, FL}
\email{mcdonald@ncf.edu}

\begin{document}

\title[Torsional rigidity and isospectrality]{Torsional rigidity, isospectrality and quantum graphs}

\begin{abstract}
We study torsional rigidity for graph and quantum graph analogs
of well-known pairs of isospectral non-isometric planar domains.  We
prove that such isospectral pairs are distinguished by torsional
rigidity.  
\end{abstract}
\maketitle

\section{Introduction}

In 1966, Marc Kac offered a poetic formulation of what came to be a
much studied problem involving the geometry of planar domains:  Does
Dirichlet spectrum determine a planar domain up to isometry?  Kac's
problem inspired a great deal of work much of which centered on
developments involving heat invariants.  The appearance of Sunada's
work involving the construction of isospectral metrics for Riemannian
manifolds provided powerful new methods for investigating diverse
spectral phenomena; in particular, Sunada's method led to the
discovery of a number of negative results for inverse spectral
problems.  In 1992 Gordon, Web and Wolpert used an extension of
Sunada's method to construct a pair of isospectral, non-isometric
planar domains \cite{GWW}.  Soon after, Buser, Conway Doyle and
Semmler produced a collection of families of pairs of isospectral
non-isometric planar polygonal domains (we refer to these
examples as {\it BCDS pairs}).  Their construction, again going back to
Sunada, involves reflecting a ``seed'' triangle across edges (see
\cite{BCDS}). An example of their construction can  
be found in section 3 below and plays an important role in our work.  

Given these negative results in reference to Kac's
question,\footnote{For smoothly bounded domains Kac's problem is
  still open.} one might ask: What are good invariants for planar
domains and/or quantum graphs? 

One approach to finding new invariants is to use data which, while not
spectral, arise from the normalized Dirichlet eigenfunctions in a natural
fashion.  An example is afforded by the sequence generated by counting
nodal domains.  This approach has been investigated in a recent series
of papers in which quantum graph analogs of BCDS pairs are produced
and their corresponding nodal domain counts are compared (see
\cite{BPB}, \cite{BSS}).  

Another approach to producing non-spectral invariants is given by
integration against the heat kernel.  More precisely, suppose $p_D(t,x,y)$
is the Dirichlet heat kernel for the domain $D$ and $u(t,x) = \int_D
p_D(t,x,y)dy$ is the solution of the inital value problem
\begin{align*}
\frac{\partial u}{\partial t}& =  \Delta u  \hbox{ on }
(0,\infty)\times D \\
\lim_{t\to 0}u(t,x) & =  1.
\end{align*}
Then one can define the {\it heat content of $D$}
\[
q_D(t) = \int_D u(t,x)dx
\]
as well as a sequence of invariants
\begin{equation}\label{ak1}
A_k   =  k \int_0^\infty t^{k-1} q_D(t) dt .
\end{equation}
The invariant $A_1,$ sometimes referred to as the {\it torsional 
  rigidity of the domain $D,$} can be computed via the solution of a
Poisson problem.  More precisely, if $u_1$ solves 
\begin{align*}
  \Delta u_1& = -1 \hbox{ on } D \\
 u_1 & =  0 \hbox{ on } \partial D
\end{align*}
then 
\[
A_1 = \int_D u_1(x) dx.
\]
Torsional rigidity arises in the study of elasticity and has a long
history (a source for early work is \cite{Po1}).  In addition to its
role in the mechanics of solid bodies, torsional rigidity is closely
related to the expected exit time of Brownian motion from the domain
$D$ and higher values of $k$ in (\ref{ak1}) correspond to higher moments of the exit
time (cf \cite{Mc1}, \cite{KMM} and references within).  For this reason
the sequence $\{A_k\}$ defined by (\ref{ak1}) has been labeled the {\it $L^1$-moment
  spectrum} associated to $D.$  

It is easy to see that both heat content and the moment spectrum
are invariants of the metric.  Recent work suggest they play a
valuable role in understanding fundamental geometric properties of a
given ambient space (for example, isoperimetry \cite{HMP1}
\cite{HMP2}, \cite{Mc2}), and one might reasonably ask how such
objects compare to Dirichlet spectrum as tools to classify domain
behavior up to isometry.  Such a study has been initiated in
\cite{MM1}.  For the purpose of this note, the main results of \cite{MM1} can be summarized as follows:  
\begin{enumerate}
\item For smoothly bounded domains, the moment spectrum determines the
  heat content. 
\item For generic smoothly bounded domains, the moment spectrum
  determines the Dirichlet spectrum.
\end{enumerate}

It is a result of Gilkey that heat content cannot distinguish
isospectral planar domains constructed via the Sunada method
\cite{Gi1}.  On the other hand, in \cite{BDK}, the authors establish that heat content distinguishes the isosspectral domains constructed by Chapman and in \cite{MM2} the authors construct
combinatorial analogs of BCDS pairs and explicit check that their heat
contents differ (at the fifth coefficient).  The main result of this
paper is that with respect to heat content, quantum graph analogs of
BCDS pairs behave more like combinatorial graphs than planar domains, and that torsional rigidity is sufficient to distinguish well-known examples of isospectral pairs.  More precisely, we have 

\begin{thrm}\label{theorem1}Let $D_1$ and $D_2$ be the $7_1$-pair
  of isopectral non-isometric planar polygonal domains as constructed
  by Buser, Conway, Doyle and Semmler and let $G_1$   and $G_2$ be their
  quantum graph analogs (cf section 3 below).  Then $G_1$ and $G_2$
  are isospectral, non-isometric and distinguished by their
  torsional rigidity.    
\end{thrm}

In addition to establishing Theorem 1, we prove 

\begin{thrm}\label{theorem2} The moment 
spectrum of a quantum graph with Dirichlet standard boundary conditions (cf section 2 below) determines the heat content of the quantum graph.
\end{thrm}
As an immediate corollary, we conclude that our isospectral
non-isometric pairs are distinguished by their heat content. Our
technique applies to the other families of isospectral pairs
constructed in \cite{BCDS}.    

The remainder of this note is organized as follows.  In the second 
section we establish the required background information concerning
quantum graphs, establish our notational conventions and prove Theorem
\ref{theorem2}.  In the third section we review the Sunada
construction in the context of planar domains following Buser, Conway,
Doyle and Semmler and include the construction of isospectral
non-isometric quantum graph analogs of BCDS pairs.  In the fourth
section of the paper we explicitly compute the torsional rigidity for
a pair of isospectral non-isometric quantum graphs and check they
differ.  We end the paper with a previously studied example
(\cite{MM2}) of a pair of  isospectral non-isomorphic weighted
combinatorial graphs which arise as analogs of BCDS pairs and check
that they too are distinguished by their torsional rigidity. 

\section{Quantum graphs}

Let $G$ be a graph with finite vertex set $V$ and edge set $E.$  For $v\in V,$ denote by $d(v)$ the cardinality of the
set $E_v = \{e\in E : v\in e\}.$  We will often find it convenient
to assign an orientation to edges; that is, for each edge
$e=\{u,v\},$ choose ${\mathcal O}(e)$ to be either $u$ or $v.$

By a {\it path} in $G$ we will mean a sequence of vertices,
$v_{i_1}, v_{i_2},\dots v_{i_m},$ with $\{v_{i_j},v_{i_{j+1}}\} \in
E.$  We say $G$ is connected if every pair of vertices can be
connected by a path in $G.$  We restrict our attention to graphs
which are connected and contain no edges of the form $\{v,v\}.$

We endow $G$ with a metric structure as follows.  For each edge
$e\in E,$ choose a positive real number $l_e.$  Identify each edge
$e\in E$ with the interval $[0,l_e],$ where $0$ is identified with
${\mathcal O}(e)$ and $l_e$ is identified with the remaining vertex
of $e.$  This identification allows us to introduce a natural
coordinate on each edge; the coordinate $x_e$ along the interval
$[0,l_e].$  These coordinates give rise to a natural metric
structure for $G.$  The pair $(G,\{l_e\}_{e\in E})$ is called a {\it
metric graph}.  Because $G$ and $l_e$ are finite, the resulting
metric graph is compact.

We identify functions on $(G,\{l_e\}_{e\in E})$ with functions along
the open edges together with values at each vertex:  For $\phi:G \to
{\mathbb C},$ we write $\phi = \oplus_{e\in E} \phi_e.$  
The metric structure on each edge gives rise to a natural Hilbert
space associated to $(G,\{l_e\}_{e\in E}).$  Write
\[
{\mathcal H}_e = L^2([0,l_e])
\]
where the inner product associated to edge $e$ is defined by 
\[
\langle f, g\rangle_e = \int_0^{l_e} f(x_e)\overline{g(x_e)} dx_e.
\]
Let 
\[
{\mathcal H} = \oplus_{e\in E} {\mathcal H}_e
\]
and denote the inner product on ${\mathcal H}$ by 
\begin{align*}
\langle f, g\rangle & =  \int_G f(x)\overline{g(x)} dx \\
  & \equiv  \sum_{e\in E} \langle f, g\rangle_e.
\end{align*}
In particular, 
\[
\langle 1, 1\rangle = L(G)
\]
where $L(G)$ is the total length of the graph $G.$

There is a
natural differential operator acting on functions on the interior of
each edge: $\frac{d^2}{dx_e^2}.$  We wish to extend this operator to
a self-adjoint operator on $L^2$-functions on the metric graph.  We sketch
the approach developed in \cite{KPS1} (see \cite{KPS1} for details).  

For each $e\in E,$ let ${\mathcal D}_e$ denote
\[
{\mathcal D}_e = \{\phi_e\in  {\mathcal H}_e : \phi_e, \
\phi^\prime_e \hbox{ are absolutely continuous, } \phi^{\prime
\prime}_e \in {\mathcal H}_e\}.
\]
Let
\[
{\mathcal D}_e^0 = \{\phi_e\in  {\mathcal D}_e :
\phi_e(0)=\phi_e(l_e)= \phi_e^\prime(0)= \phi^\prime_e (l_e) =0 \}
\]
and set
\begin{align*}
{\mathcal D} & =  \oplus_{e\in E} {\mathcal D}_e \\
{\mathcal D}^0 & =  \oplus_{e\in E} {\mathcal D}^0_e.
\end{align*}
Then the operator $\Delta^0$ which acts on $\phi \in {\mathcal D}^0$
according to
\[
(\Delta^0\phi)_e (x) = \frac{d^2 \phi_e}{dx_e^2}(x_e)
\]
is closed, symmetric and densely defined.  We seek to impose
boundary conditions which give self-adjoint extensions of $\Delta^0$
to ${\mathcal H}.$

We begin by noting that boundary values of functions along edges lie
in the space ${\mathcal B} = {\mathbb C}^{|E|} \times {\mathbb
C}^{|E|} $ where the order of the components is fixed by the
orientation.  More precisely, given $\phi \in {\mathcal D},$ let
$\underline{\phi} \in {\mathcal B}$ be given by
\[
\underline{\phi} = (\{\phi_e(0)\}_{e\in E}, \{\phi_e(l_e)\}_{e\in
E}).
\]
Boundary conditions can then be described by a pair of linear
operators $A, B: {\mathcal B} \to {\mathcal B}$ satisfying
\[
A\underline{\phi} + B {\underline \phi^\prime} = 0.
\]
Those pairs of linear operators $(A,B)$ leading to self-adjoint
extensions of $\Delta^0$ have been classified by Kostrykin, Potthoff
and Schrader \cite{KPS1} using a symplectic formalism going back to
at least to Novikov.  We are interested in a single special case, {\it
Dirichlet standard boundary conditions,} which we now describe.

\begin{definition}\label{sbc} Denote by DSBC the collection of
  continuous functions defined by 
\begin{align}
\hbox{DSBC} & =  \left\{ \phi \in {\mathcal D}: \begin{array}{ll}
                                       \phi(v) = 0 & \hbox{ if }
                                       d(v) =1 \\
                                       \sum_{e\in E_v} 
                                       \frac{d\phi}{d\nu_e}(v) = 0 &
                                       \hbox{ if } d(v)\neq 1. \end{array}
                       \right\}\label{Standard}
\end{align}
where the derivative occurring in (\ref{Standard}) is always
directed into the vertex $v$ and, as before, $E_v = \{e\in E: v\in
e\}$ is the collection of edges on which $v$ is incident.
\end{definition}

We note that the literature often refers to "standard boundary
conditions" as those involving Kirchoff conditions at internal
vertices (ie vertices of degree at least 2) and Neumann conditions
at boundary vertices (ie vertices of degree 1).

We can describe standard boundary conditions using pairs of matrices
$(A,B)$ as follows:  For each vertex $v$ with $d(v)>1,$ introduce
$d(v)\times d(v)$ matrices $A$ and $B$ (cf also \cite{KPS1}, \cite{BPB})
\begin{align*}
A(v) = \left( \begin{array}{cccccc}
         1 & -1& 0 & \cdots & 0 & 0 \\
         0 & 1 & -1 & \cdots & 0 & 0 \\
         0 & 0 & 1 & \cdots & 0 & 0 \\
         \vdots & \vdots & \vdots & \cdots & \vdots & \vdots \\
         0 & 0 & 0 & \cdots & 1 & -1 \\
         0 & 0 & 0 & \cdots & 0 & 0 \end{array} \right) &  
B(v) = \left( \begin{array}{cccccc}
         0& 0& 0 & \cdots & 0 & 0 \\
         0 & 0 & 0 & \cdots & 0 & 0 \\
         0 & 0 & 0 & \cdots & 0 & 0 \\
         \vdots & \vdots & \vdots & \cdots & \vdots & \vdots \\
         0 & 0 & 0 & \cdots & 0 & 0 \\
         1 & 1 & 1 & \cdots & 1 & 1 \end{array} \right).
\end{align*}
For vertices of degree 1, define $1\times 1$ matrices $A= 1$ and
$B=0.$  The matrices $A$ and $B$ referenced above then take a block
form:  
\begin{align*}
A = \oplus_{v\in V} A(v) &  \hspace{1in} &  = \oplus_{v\in V} B(v)
\end{align*}
We note that Dirichlet standard boundary conditions are local in the
sense of \cite{KPS3} and do not mix derivatives with function values.  

We note that Dirichlet standard boundary conditions force continuity
across all internal vertices for all elements in the corresponding
domain of the self-adjoint extension of $\Delta^0.$  Thus, it is the
choice of boundary conditions which ultimately ties the geometry and
topology of the quantum graph to the analysis of the corresponding
Laplace operator.

We will denote the self-adjoint extension of $\Delta^0$ with
Dirichlet standard boundary conditions as $\Delta.$   The pair consisting of the metric graph $G$ and the operator $\Delta$ is the quantum graph we study.  

The spectrum of $\Delta$
is discrete and of finite multiplicity.  We will denote the spectrum
of $\Delta$ by $\hbox{spec}(\Delta).$  We will assume there is at
least one boundary vertex and we will list elements of
$\hbox{spec}(\Delta)$  in increasing order with multiplicity:
\[
 0 < \lambda_1 < \lambda_2 \leq \lambda_3 \leq \dots
\]
There is a great deal known about the spectrum for arbitrary
self-adjoint boundary conditions.  For our purposes it is important
that there is a heat
kernel which can be written as 
\begin{equation}\label{heatk1}
p_G(t,x,y) = \sum_{\lambda \in \hbox{spec}(\Delta)} e^{-\lambda
  t}\phi_\lambda(x)\phi_\lambda(y) .
\end{equation}
The heat content of $G$ can then be defined via integration as for
domains: 
\begin{align}
q_G(t) & =  \int_G\int_G p_G(t,x,y) dx dy\nonumber \\
 & =  \sum_{\lambda \in \hbox{spec}(\Delta)} e^{-\lambda
  t}\left(\int_G \phi_\lambda(x)dx\right)^2. \label{heatc1}
\end{align}
We let $\hbox{spec}^*(\Delta)$ denote the set of values defined by
the Dirichlet standard boundary conditions (ie we disregard
multiplicity) for which the corresponding eigenspace is not orthogonal to constant functions.  Then we can write the heat content in a fashion which
makes explicit its non-spectral nature:
\begin{equation}\label{heatc2}
q_G(t)   =  \sum_{\lambda \in \hbox{spec}^*(\Delta)} a_\lambda^2 e^{-\lambda
  t}
\end{equation}
where $a_\lambda^2$ is the square of the $L^2$-norm of the projection
of the constant function 1 on the eigenspace defined by
$\lambda.$  

The moments of $q_G(t)$ define a sequence associated to the quantum
graph $G.$  More precisely, as in the introduction, for $k$ a natural number, set  
\begin{equation*}
A_k = k\int_0^\infty t^{k-1} q_G(t) dt.
\end{equation*}

The moment spectrum of the graph $G$ is the sequence
$\{A_k\}_{k=1}^\infty.$  

From the definition it is clear that the heat content determines the
moment spectrum.  The converse (Theorem \ref{theorem2} above) is also true:

\begin{thrm} The moment spectrum of the quantum graph $G$
  determines the heat content of the quantum graph $G.$
\end{thrm}

{\sc Proof}  The proof follows the case of domains, where the same
result holds (see \cite{MM1}, \cite{Mc2}).  For completeness we
provide the argument given in \cite{MM1}.

It suffices to show that the moment spectrum determines both
$\hbox{spec}^*(\Delta)$ and the collection of coefficients
$\{a_\lambda^2\}.$  The important observation is that the moment spectrum describes special values of the Mellin transform of the heat content.  More precisely, using (\ref{heatc2}), the Mellin transform of the heat content is given by the Dirichlet series
\begin{equation}\label{zeta1}
\zeta_G(s) = \sum_{\lambda \in \hbox{spec}^*(\Delta)} a_\lambda^2 \left(\frac{1}{\lambda}\right)^s
\end{equation}
which converges for real part of $s$ nonnegative and admits a meromorphic extension to the plane, with poles at the negative half-integers (see \cite{MM1}).  The moment spectrum then satisfies
\begin{equation}\label{zeta2}
\zeta_G(n) = \frac{A_n}{\Gamma(n+1)}
\end{equation}

Using (\ref{zeta1}) and (\ref{zeta2}), there is a recursion:  Write
$\hbox{spec}^*(\Delta) = \{\mu_n\}_{n=1}^\infty $ with $\mu_n$ strictly
increasing.  Then
\[
\mu_1 = \sup \left\{ \mu\geq 0: \limsup_{n\to \infty} (\mu)^n
\frac{A_n}{\Gamma(n+1)} < \infty\right\}
\]
and
\[
a^2_{\mu_1} = \limsup_{n\to \infty} (\mu_1)^n
\frac{A_n}{\Gamma(n+1)}.
\]
Having determined $\mu_j$ and $a_{\mu_j}^2$ for $j<k,$ we have
\[
\mu_k = \sup \left\{ \mu\geq 0: \limsup_{n\to \infty} (\mu)^n
\left(\frac{A_n}{\Gamma(n+1)} -\sum_{j=1}^{k-1} a_{\mu_j}^2
\left(\frac{1}{\mu_j}\right)^n\right) < \infty\right\}
\]
and
\[
a^2_{\mu_k} = \limsup_{n\to \infty}
(\mu_k)^n\left(\frac{A_n}{\Gamma(n+1)} -\sum_{j=1}^{k-1}
a_{\mu_j}^2 \left(\frac{1}{\mu_j}\right)^n\right).
\]
This shows that both $\hbox{spec}^*(\Delta)$ and the length partition of
$G$ are determined by the $L^1$-moment spectrum.  Using the spectral
representation of the heat content (\ref{heatc1}), it is clear that
$\hbox{spec}^*(\Delta)$ and the coefficients $a_\lambda^2$ determine the heat content.
This proves that the $L^1$-moment spectrum determines heat content and
completes the proof of Theorem \ref{theorem2}.
 
\endproof

Our primary interest is in $A_1,$ which, as in the case of smooth
domains, can be described via a solution of a Poisson problem on the
graph $G.$  More precisely, Let $u_1$ solve 
\begin{align}
\Delta u_1 & =  -1 \hbox{ on } {\mathrm{interior}}(G) \nonumber \\
u_1 & =  0 \hbox{ on } \partial G. \label{poisson}
\end{align}
We can apply Fubini's theorem to interchange the order of
integration in the definition of the heat content and $u_1(x)$ to obtain
\begin{equation}\label{tt2}
A_1 = \int_G u_1(x) dx.
\end{equation}
Given any quantum graph with Dirichlet standard boundary conditions, 
the computation of the right-hand-side of (\ref{tt2}) is an
exercise in linear algebra.  We compute torsional rigidity for certain
pairs of quantum graphs constructed in the next section.  

\section{The Sunada construction, BCDS pairs and isospectral quantum graphs}

As mentioned in the introduction, Buser, Conway, Doyle and Semmler
have constructed families of pairs of isospectral planar domains
that are not isometric.  Their construction and investigation of the
resulting pairs is spectacularly simple.  To construct families of
pairs of isospectral planar domains,
\begin{itemize}
\item  Fix a seed triangle with edges labeled;
\item  Reflect the seed about its edges to produce a second
generation of labeled triangles;
\item  For each of the progeny, reflect and label across edges
determined by two rules derived from subgroups which are "nearly
conjugate" in Sunada's sense
\item  Iterate
\item  Interpret the edges as distinct lengths of an acute triangle. 
\end{itemize}
An example of the construction is given in Figure 1: the
so-called $7_1$ pairs (for complete details see \cite{BCDS}).

\begin{figure}\label{fig1}
\begin{center}
\includegraphics{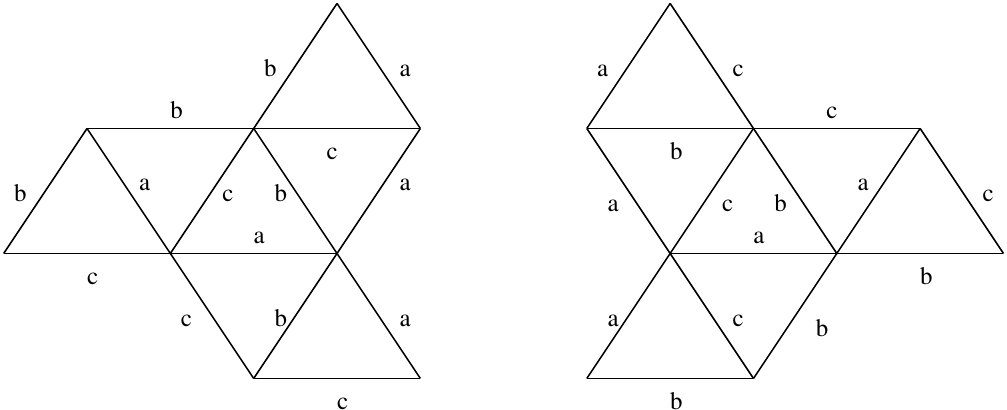}
\end{center}
\caption{Pairs of isospectral non-isometric planar domains}
\end{figure}

Having constructed families of pairs as above, Buser et al prove
that the associated pairs are isospectral using the technique of
eigenfunction transplantation. Given an eigenfunction $f$ on one
domain, say $D_1,$ they use the process by which the domain is
constructed to naturally partition the domain into a collection of
isometric subdomains, say $T_{1,i}$ where $i$ indexes the number of
triangles occurring in the construction, say $1\leq i \leq n.$  The
same procedure carried out on the domain pair, say $D_2,$ results in a
collection $T_{2,i}, \ 1\leq i \leq n,$ of isometric subdomains of
$D_2.$  Having partitioned both domains into collections of
subdomains, they restrict the given eigenfunction to each subdomain to
produce an associated collection of functions, $f_i, \ 1\leq i \leq
n,$ on isometric copies of the same subdomain.  They then use the
geometry of the underlying domains $D_1$ and $D_2$ (ie, how 
the domains are constructed from the collections of subdomains) to
construct a linear map which prescribes how to build an
associated eigenfunction on the domain $D_2$ by combining the function
elements $f_i$ on each subdomain of $D_2.$  The process by which an
explicit transplantation map is constructed is described neatly in the
references \cite{BCDS} and \cite{Ch}.  For our purposes it is
important to note that the process is essentially combinatorial in the
sense that it depends only on the following: 
\begin{itemize}
\item  The eigenfunction $f$ vanishes at the boundary of $D_1.$  
\item  To smoothly continue the eigenfunction $f$ across the boundary
  of the domain $D_1,$ the reflection principle requires that the
  value of $f(x^*)$ be prescribed to be $-f(x)$ where $x$ is the
  reflection image of $x^*.$ 
\item  If $T_{1,i_1}$ and $T_{1,i_2}$ are subdomains which share a
  boundary internal to $D_1,$ the associated eigenfunction elements
  $f_{i_1}$ and $f_{i_2}$ and their respective normal derivatives must
  match along the common boundary.
\end{itemize}

There are also techniques to produce {\it
  quantum graph} analogs \cite{BPB}, \cite{SS}, \cite{BSS}. To
proceed, we recall the required definitions.  

Following earlier work (cf \cite{SS}, \cite{BSS},
\cite{BPB} and references therein) Band, Parzanchevski and Ben-Shach
have described an extension of eigenfunction transplantation in the
context of quantum graphs \cite{BPB}.  They provide a method which
associates to every family of pairs constructed in \cite{BCDS}, a family of
pairs of quantum graphs which are isospectral but not isometric. The
construction technique is straightforward, as is the check that the
objects satisfy the required conditions.  We sketch the process
below and produce the pairs associated to the $7_1$ pairs in Figure 2.

For each pair of isospectral nonisometric domains found in
\cite{BCDS},
\begin{itemize}
\item  For each triangle appearing in the BCDS construction, introduce
  a 3-star graph consisting of three edges, with edges labeled with a
  label from the edges of the corresponding triangle, joined at a
  central vertex. 
\item  Join edges in 3-stars at a common degree two vertex if
corresponding edges in the triangles to which they are associated
overlap.
\end{itemize}
We carry out the process on quantum graphs corresponding to $7_1$
pairs in Figure \ref{fig3}.  Note that for the graphs to be
isometric, it would be necessary to map the 3-stars corresponding to
generating triangles onto each other. 

\begin{figure}\label{fig2}
\begin{center}
\includegraphics{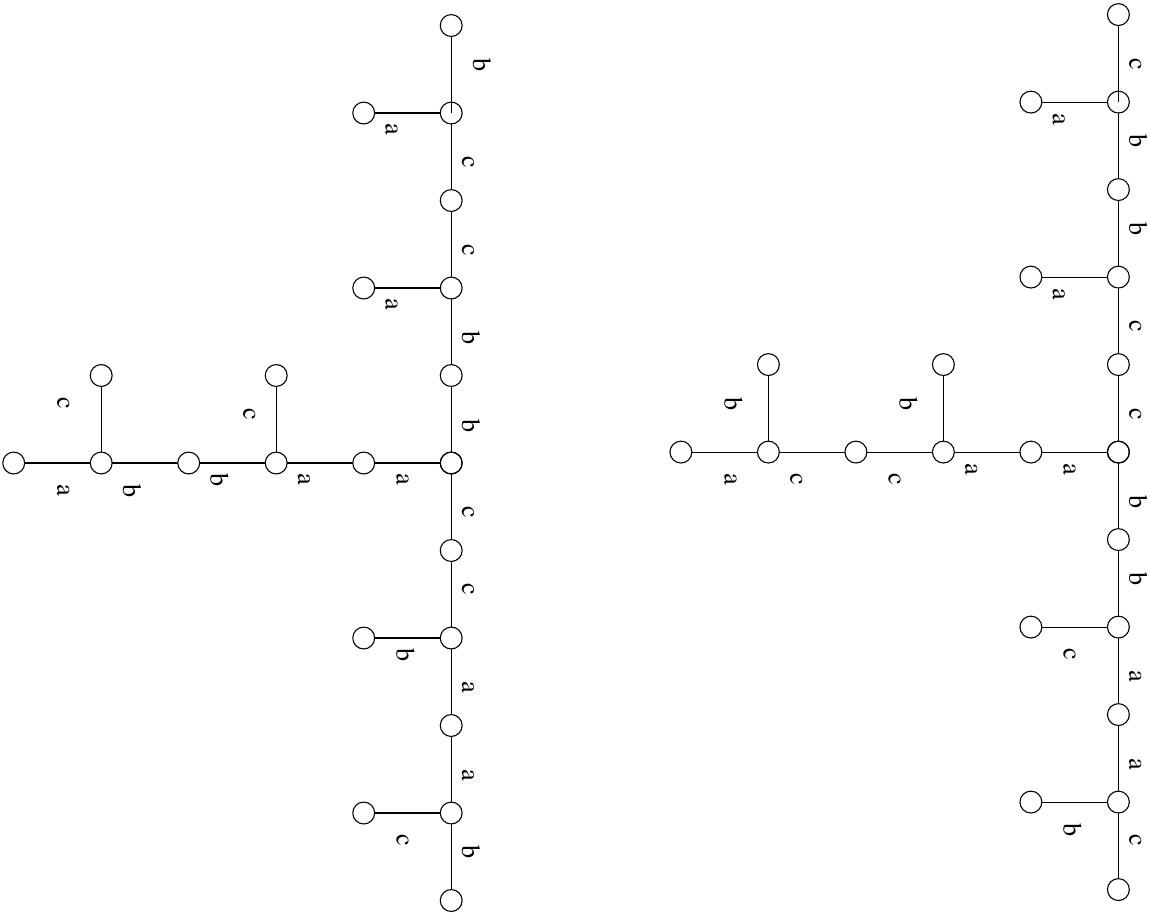}
\end{center}
\caption{Quantum graph analogs of pairs of isospectral non-isometric
  planar domains}  
\end{figure}

For Neumann standard boundary conditions, the quantum graph analogs
are shown to be isospectral and non-isometric in \cite{BSS} (see also
\cite{BPB}).  This is accomplished by showing that the graphs share
the same secular equation.  The same argument can be made to work in
the current context, but instead we provide an explicit
transplantation map.  In fact, it's easy to check that the transplantation
map used to establish that the BCDS pairs are isometric works to show
that the associated quantum graphs pairs are isometric.  We do this as
follows:
\begin{itemize}
\item  Label the edges of the quantum graph $G_1$ and the corresponding
partner $G_2.$
\item  Given an eigenfunction on $G_1,$ restrict to 3-star subgraphs
  to obtain a collection of functions defined on isometric 3-stars.  
\item  Use the linear map constructed following \cite{BCDS} to
  construct the required transplantation map taking eigenfunctions of
  $G_1$ to eigenfunctions of $G_2.$   
\end{itemize}
We can formalize these observations with 
\begin{prop}\label{71pairs}The quantum graph (with Dirichlet
  standard boundary conditions) analogs of the 
 pairs constructed in \cite{BCDS} are isospectral and non-isometric.
\end{prop}
{\it Proof}  By construction, for each domain constructed in
\cite{BCDS} the corresponding quantum graph is obtained by replacing
subdomains with 3-stars, gluing 3-stars along boundaries which
correspond to boundaries which are glued in subdomains.  The
transplantation map for domains may be represented by a matrix acting
on function elements defined on subdomains, all of which can be
identified.  This defines a linear map from eigenfunctions on $G_1$ to
functions on $G_2,$ which we will write as ${\mathcal L}.$  Give an eigenfunction
$\phi$ on $G_1$ corresponding to eigenvalue $\lambda,$ to see that the
image ${\mathcal L}\phi$ is an eigenfunction on $G_2$ first note that on the
interior of every edge, $\Delta {\mathcal L}\phi = \lambda {\mathcal L}\phi$ by linearity.
Thus, it suffices to check that ${\mathcal L}\phi$ has value zero at the boundary
of $G_2$ and that ${\mathcal L}\phi$ behaves as it should across glued edges.
But for this to be the case requires only the three properties listed
in our description of the transplantation process given above.  Since
the gluing relations for quantum graphs are precisely those that are
prescribed by the domains and the boundary conditions are Dirichlet in
both cases, these three properties must hold for the image of $\phi$
under ${\mathcal L}.$  In particular, ${\mathcal L}\phi$ must be an eigenfunction
corresponding to eigenvalue $\lambda.$  

We give the required transplantation map in Figure 3.  

\begin{figure}\label{fig3}
\begin{center}
\includegraphics{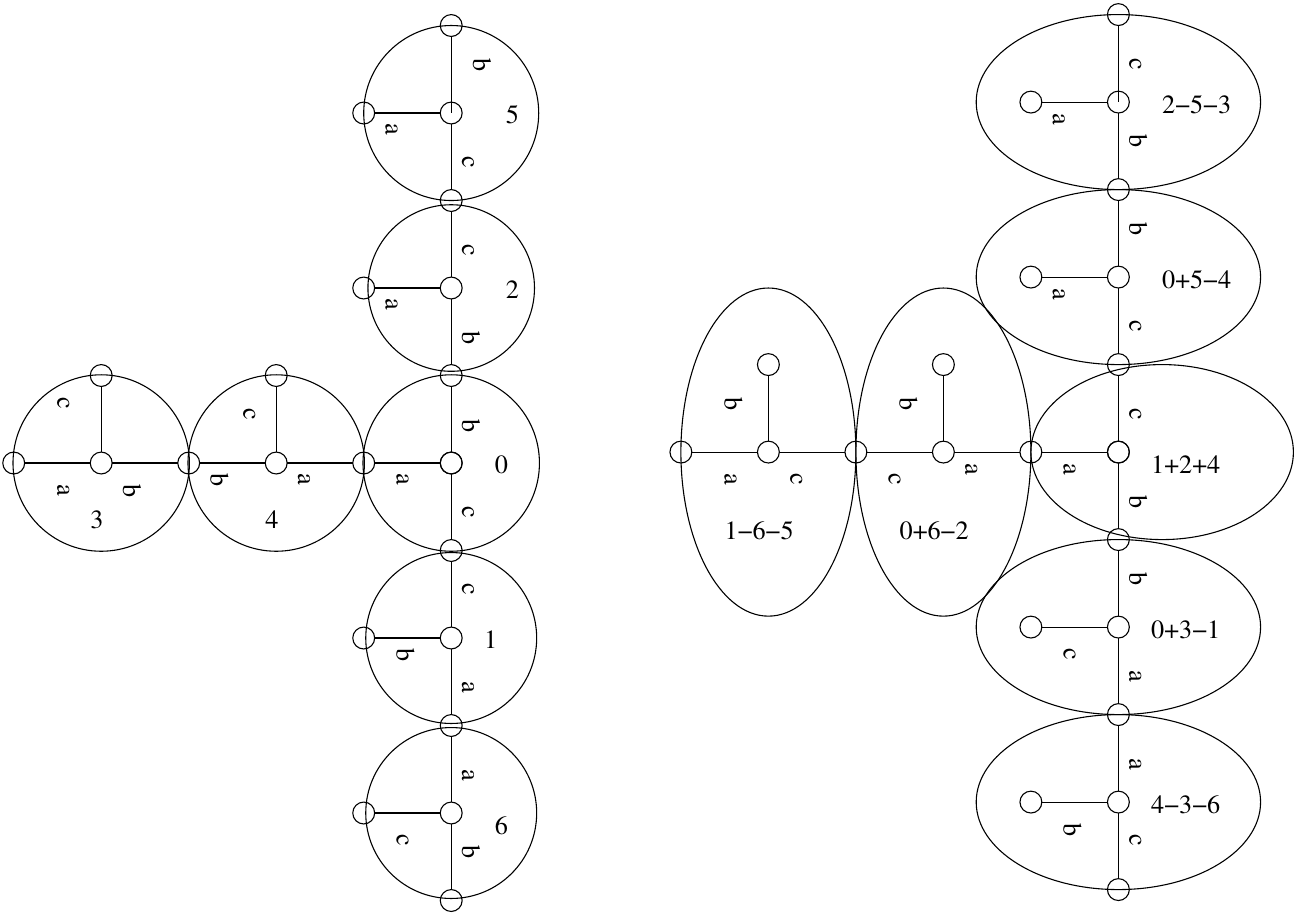}
\end{center}
\caption{Eigenfunction transplantation for quantum graphs}
\end{figure}

\section{Computing Torsional Rigidity}

We are now in a position to prove Theorem \ref{theorem1} for the $7_1$
pairs of Figure 2.  We begin by noting that we can
eliminate vertices of degree 2 and replace the corresponding two
edges by a single edge of length equal to the sum of the
corresponding edges (following Friedlander \cite{Fr1}, we call such quantum
graphs {\it clean}).  This results in the pair of quantum graphs
given in Figure 4.  Note that there is a symmetry: interchanging $b$ and $c$ in the two graphs maps one graph onto the other.  We use this symmetry in the computation of the torsional rigidity of
each graph.

Suppose we denote by $x_{i,j}$ for the coordinate on the $j$th edge of the $i$th graph.  Then we must solve
\begin{align*}
\frac{d^2 u_{i,j}}{dx_{i,j}^2} & =  -1 
\end{align*}
with boundary conditions determined by our choice of Dirichlet standard boundary conditions.  The solution on each edge is of the form
\begin{equation}\label{torsion2}
u_{i,j}(x_{i,j}) = -\frac{1}{2} x_{i,j}^2 +\alpha_{i,j} x_{i,j} +
\beta_{i,j}
\end{equation}
where the constants $\alpha_{i,j}$ and $\beta_{i,j}$ are to be
determined from the Dirichlet standard boundary conditions.  Thus,
for each graph there are 30 constants to be determined using nine
Dirichlet boundary conditions at vertices of degree 1, seven
Kirchoff boundary conditions at vertices of degree 3, and fourteen
continuity conditions at vertices of degree 3.  This determines a
linear system $L{\bf a} ={\bf v}$ where the coefficient matrix
$L$ is a sparse $21\times 21$ matrix and the constant vector ${\bf v}$
is defined in terms of continuity constraints and conservation at
internal nodes.  The structure of the matrix will be determined by our choice of edge labels and orientation.  We choose the first nine edges to be those incident to boundary vertices and we orient the edges to be inward facing relative to the boundary.  With this choice the coefficients $\beta_{i,j} = 0 $ for $i= 1, \ 2$ and $1 \leq j\leq 9.$  With the remaining  edges oriented to point toward the "central node", the representation of $L$ for $G_1$ as pictured in figure
4, is given by the coefficient matrix $L_1$ 
\tiny
\[ \left(\begin{array}{ccccccccccccccccccccc}
 c & 0& 0& 0& 0& 0& 0& 0& 0& 0& 0& 0& 0& 0& 0& -1& 0& 0& 0& 0& 0 \\
 0 & a& 0& 0& 0& 0& 0& 0& 0& 0& 0& 0& 0& 0& 0& -1& 0& 0& 0& 0& 0 \\
 0 & 0& b& 0& 0& 0& 0& 0& 0& 0& 0& 0& 0& 0& 0& 0& 0& 0& 0& -1& 0 \\
 0 & 0& 0& a& 0& 0& 0& 0& 0& 0& 0& 0& 0& 0& 0& 0& 0& 0& 0& -1& 0 \\
 0 & 0& 0& 0& b& 0& 0& 0& 0& 0& 0& 0& 0& 0& 0& 0& 0& 0& 0& 0& -1 \\
 0 & 0& 0& 0& 0& c& 0& 0& 0& 0& 0& 0& 0& 0& 0& 0& 0& 0& 0& 0& -1 \\
 0 & 0& 0& 0& 0& 0& a& 0& 0& 0& 0& 0& 0& 0& 0& 0& -1& 0& 0& 0& 0 \\
 0 & 0& 0& 0& 0& 0& 0& b& 0& 0& 0& 0& 0& 0& 0& 0& 0& 0& -1& 0& 0 \\
 0 & 0& 0& 0& 0& 0& 0& 0& c& 0& 0& 0& 0& 0& 0& 0& 0& -1& 0& 0& 0 \\
 0 & 0& 0& 0& 0& 0& 0& 0& 0& 2b& 0& 0& 0& 0& 0& 1& -1& 0& 0& 0& 0 \\
 0 & 0& 0& 0& 0& 0& 0& 0& 0& 0& 2c& 0& -2a& 0& 0& 0& 1& 0& -1& 0& 0 \\
 0 & 0& 0& 0& 0& 0& 0& 0& 0& 0& 0& -2b& 2a& 0& 0& 0& 0& -1& 1& 0& 0 \\
 0 & 0& 0& 0& 0& 0& 0& 0& 0& 0& 0& 0& 0& 2c& 0& 0& 0& 0& -1& 1& 0 \\
 0 & 0& 0& 0& 0& 0& 0& 0& 0& 0& 0& 0& 0& 0& 2a& b& 0& -1& 0& 0& 1 \\
 1 & 1& 0& 0& 0& 0& 0& 0& 0& -1& 0& 0& 0& 0& 0& 0& 0& 0& 0& 0& 0 \\
 0 & 0& 1& 1& 0& 0& 0& 0& 0& 0& 0& 0& 0& -1& 0& 0& 0& 0& 0& 0& 0 \\
 0 & 0& 0& 0& 1& 1& 0& 0& 0& 0& 0& 0& 0& 0& -1& 0& 0& 0& 0& 0& 0 \\
 0 & 0& 0& 0& 0& 0& 1& 0& 0& 1& -1& 0& 0& 0& 0& 0& 0& 0& 0& 0& 0 \\
 0 & 0& 0& 0& 0& 0& 0& 1& 0& 0& 0& 0& -1& 1& 0& 0& 0& 0& 0& 0& 0 \\   
 0 & 0& 0& 0& 0& 0& 0& 0& 1& 0& 0& -1& 0& 0& 1& 0& 0& 0& 0& 0& 0 \\
 0 & 0& 0& 0& 0& 0& 0& 0& 0& 0& 1& 1& 1& 0& 0& 0& 0& 0& 0& 0& 0 
  \end{array} \right)
\]
\normalsize
The continuity and conservation constraints define a vector, ${\bf
  v}_1,$ given by  
\small
\begin{align*}
{\bf v}_1 & =  (c^2/2, a^2/2, b^2/2,  a^2/2, b^2/2, c^2/2, a^2/2, b^2/2, c^2/2, 2 b^2, 2c^2-2a^2, 2a^2-2b^2, 2c^2, 
  2a^2, \\
         &   a+c, a+b, b+c, a+2b, b+2c,  2a +c, 2a+2b+2c)^T
\end{align*}
\normalsize
where ``$T$'' denotes transpose.  To obtain the data for the second graph, simply interchange $b$ and $c.$

We can use these systems to carry out the required calculation
\small
\begin{align}
\int_{G_1} u_{1}(x) dx - \int_{G_2} u_2(x) dx & =  \sum_{j=1}^{15}
\left[\int_0^{l_{1,j}}u_{1,j} (x_{1,j}) dx_{1,j} -
  \int_0^{l_{2,j}}u_{2,j} (x_{2,j}) dx_{2,j}\right] \nonumber \\ 
 & =  \frac{1}{2} \sum_{j=1}^{15} \left(\alpha_{1,j} l_{1,j}^2 -  \alpha_{2,j}
 l_{2,j}^2 \right)  + \sum_{j=1}^{15}\left( \beta_{1,j}l_{1,j} -
 \beta_{2,j}l_{2,j}\right) \nonumber \\ 
  & =  \langle L_1^{-1}{\bf v}_1,{\bf l}_1\rangle - \langle L_2^{-1}{\bf
     v}_2,{\bf l}_2 \rangle \label{tor33}
\end{align}
\normalsize
where the vectors ${\bf l}_i$ are constructed from the lengths of the
edges of the graphs $G_i:$
\small
\begin{align*}
{\bf l}_1 & =  1/2(c^2, a^2, b^2, a^2, b^2, c^2, a^2, b^2, c^2, 4 b^2, 4 c^2, 4b^2, 4 a^2, 4 c^2, 4 a^2, 4 b, 4 c, 4 b, 
  4a, 4c, 4a)^T 
\end{align*}
\normalsize
and $l_2$ is obtained from $l_1$ by interchanging $b$ and $c.$  The right hand side of (\ref{tor33}) is a rational function in the variables $a,$ $b,$ and $c:$
\begin{equation}\label{sign1}
\langle L_1^{-1}{\bf v}_1, {\bf l}_1\rangle  - \langle L_2^{-1}{\bf v}_2, {\bf
    l}_2 \rangle = (a - b) (a - c) (b - c)  R(a,b,c)
\end{equation}
where the rational function $R(a,b,c)= \frac{N(a,b,c)}{D(a,b,c)}$ has a sign:
\begin{align*}
N(a,b,c) &=  -4(b c + a (b + c)) (11 a^3 (b + c) +
       b c (11 b^2 + 23 b c + 11 c^2) + \\
       &   a (b + c) (11 b^2 + 48 b c + 11 c^2) +
       a^2 (23 b^2 + 59 b c + 23 c^2))  \\
 D(a,b,c) & =  32 a^4 (b + c)^2 +
     8 b^2 c^2 (4 b^2 + 9 b c + 4 c^2) +
     8 a b c (b + c) (8 b^2 + 23 b c + 8 c^2) + \\
     &    8 a^3 (b + c) (9 b^2 + 22 b c + 9 c^2) +
     a^2 (32 b^4 + 248 b^3 c + 439 b^2 c^2 + 248 b c^3 + 32 c^4) .
\end{align*}
\normalsize 
From (\ref{sign1}) and the explicit expression of $R(a,b,c)$ it is easy to see that the difference
in torsional rigidities is nonzero whenever $a,$ $b$ and $c$ are distinct.  This proves Theorem \ref{theorem1}.

As a corollary of Theorem \ref{theorem1} and Theorem \ref{theorem2}, we obtain:

\begin{cor} The quantum graphs $G_1$ and $G_2$ are distinguished by their heat content.  
\end{cor}

\begin{figure}\label{fig4}
\begin{center}
\includegraphics{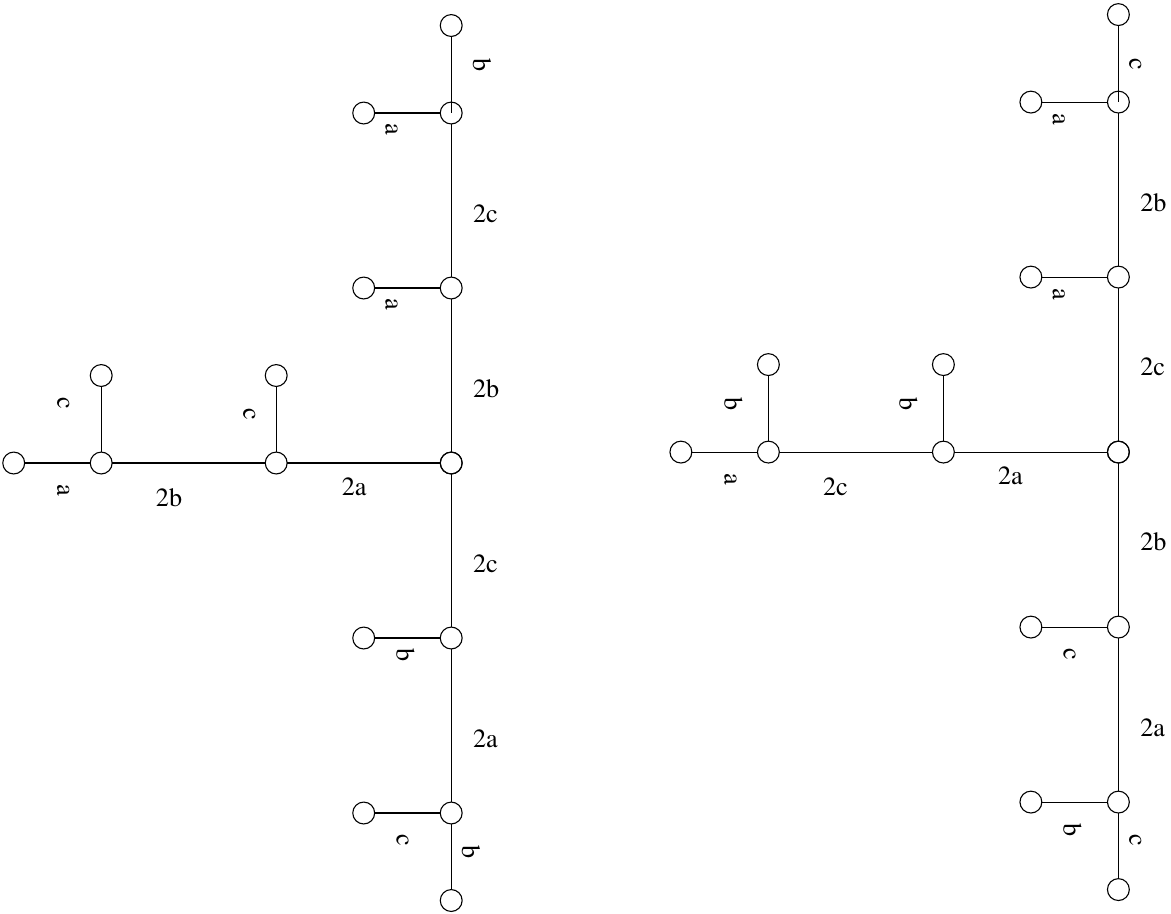}
\end{center}
\caption{Clean pairs of quantum graphs}
\end{figure}

\section{Combinatorial graphs}

In \cite{MM2} the authors construct combinatorial weighted graph
analogs of BCDS pairs.  The result of these constructions for the
$7_1$-pairs given above are the two isospectral,
non-isomorphic weighted combinatorial graphs given in Figure 5:  

\begin{figure}\label{fig5}
\begin{center}
\includegraphics[scale =0.75]{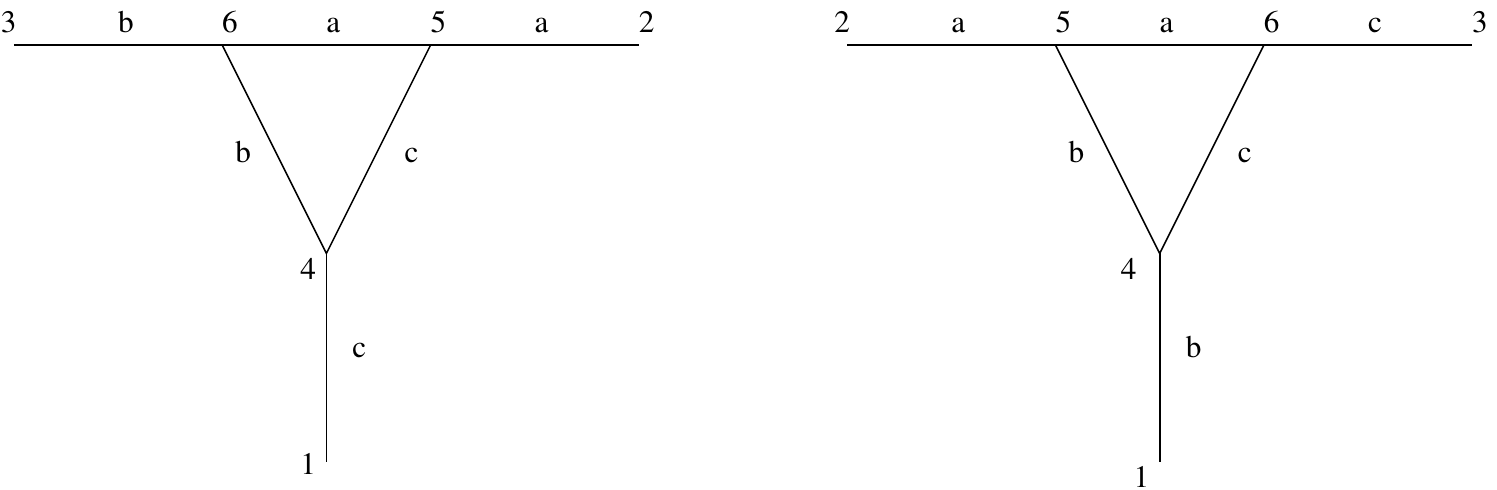}
\end{center}
\caption{Isospectral weighted combinatorial graphs}
\end{figure}

There is a natural weighted Laplacian associated to such weighted
graphs.  For functions on the vertices we write 
\[
D f(x) = \sum_{y\tilde x} W_E(x,y) f(y) - f(x)
\]
where the sum is over vertices $y$ adjacent to $x$ and $W_E(x,y)$
denotes the weight along the edge defined by the vertices $x$ and $y.$
For the above two graphs this results in Laplace operators given by 
\begin{align*}
D_1 & =  \left( \begin{array}{cccccc} 
 a + c & 0& 0& -c/2& 0& 0 \\
 0 & a + b & 0& 0& -a/2&  0 \\
 0 & 0 & b+c & 0& 0& -b/2 \\
-c/2& 0& 0& b + c& -c/2 -b/2 \\
 0& -a/2& 0 & -c/2& a + c& -a/2 \\
 0& 0& -b/2& -b/2& -a/2& a + b \end{array} \right) \\
D_2 & =  \left( \begin{array}{cccccc} 
 a + b & 0& 0& -b/2& 0& 0 \\
 0 & a + c & 0& 0& -a/2&  0 \\
 0 & 0 & b+c & 0& 0& -c/2 \\
-b/2& 0& 0& b + c& -b/2 -c/2 \\
 0& -a/2& 0 & -b/2& a + b& -a/2 \\
 0& 0& -c/2& -c/2& -a/2& a + c \end{array} \right) 
\end{align*}
It is easy to check that these matrices have the same characteristic
equation.

To compute the difference in torsional rigidity for the corresponding
combinatorial graphs we let ${\bf v} = (1,1,1,1,1,1)^T$ and we compute:
\begin{align*}
\langle D_1^{-1}{\bf v}, {\bf v}\rangle - \langle D_2^{-1}{\bf v}, {\bf v}\rangle & =
 (a-b)(a-c)(b-c) r(a,b,c)
\end{align*}
where $r(a,b,c)= \frac{n(a,b,c)}{d(a,b,c)} $ is the rational function
defined by 
\begin{align*}
n(a,b,c) & =  56( b c + a (b + c))\\
d(a,b,c) & =  (32 a^4 (b + c)^2 +  8 b^2 c^2 (4 b^2 + 9 b c + 4 c^2) + 
   8 a b c (b + c) (8 b^2 + 23 b c + 8 c^2) +\\
  &    8 a^3 (b + c) (9 b^2 + 22 b c + 9 c^2) + 
   a^2 (32 b^4 + 248 b^3 c + 439 b^2 c^2 + 248 b c^3 + 32 c^4)).
\end{align*}
In particular, for distinct values of $a,$ $b,$ and $c$ the
isospectral non-isomorphic combinatorial graphs constructed in
\cite{MM2} are distingusihed by their torsional rigidity.


\begin{thebibliography}{BCDS}

\bibitem[Br]{Br} R. Brooks (1999) {\it Non Sunada graphs.}
  Ann. Inst. Four. {\bf 49}, 707--725.  

\bibitem[BSS]{BSS} R. Band, T. Schapira and U.Smilansky (2006) {\it
  Nodal domains on isospectral quantum graphs: the resolution of
  isospectrality?} J. Phys. A  {\bf 39}  no. 45, 13999--14014 

\bibitem[BPB]{BPB} R. Band, O. Parzanchevski and G. Ben-Shach (2009)
  {\it The isospectral fruits of representation theory: quantum graphs
    and drums.} J. Phys. A {\bf 42}  no. 17, 175202, 42 pp 

\bibitem[BDK]{BDK} M. van den Berg, E. Dryden and T. Kappeler (2014) {\it Isospectrality and heat content.}
  Bull. London Math. Soc. {\bf 46} 793--808. 

\bibitem[BG]{BG} M. van den Berg and P. Gilkey (1994) {\it Heat
  content asymptotics for a Riemannian manifold with boundary.}
  J. Funct. Anal. {\bf 120} 48--71. 

\bibitem[BCDS]{BCDS} P. Buser, J. Conway, P. Doyle and K. Semmler (1994) {\it
  Some planar isospectral domains.} Inter. Math. Res. Not. {\bf 9} 391--400.

\bibitem[GWW]{GWW} C. Gordon, D. Webb and S. Wolpert (1992) {\it Isospectral
  plane domains and surfaces via Riemannian orbifolds.}
  Invent. Math. {\bf 110}, 1--22

\bibitem[Gi1]{Gi1} P. Gilkey (2009) {\it Heat content, heat trace and
  isospectrality.} New developments in Lie theory and geometry,
  Contemp. Math., {\bf 491} Amer. Math. Soc. Providence, RI 115--123. 



\bibitem[Fr1]{Fr1} L. Friedlander (2005) {\it Genericity of simple eigenvalues
  for a metric graph.}  Israel J. Math.  {\bf 146}, 149--156. 


\bibitem[KMM]{KMM} K. J. Kinateder, P. McDonald and D. Miller (1998)
  {\it Exit time moments, boundary value problems and the geometry of
    domains in Euclidean space.} Prob. Th. and Rel. {\bf 111} 469--487.

\bibitem[Mc1]{Mc1} P. McDonald (2002) {\it Isoperimetric conditions,
  Poisson problems and diffusions in Riemannian manifolds.} Potential
  Analysis {\bf 16} 115--138.

\bibitem[Mc2]{Mc2} P. McDonald (2012) {\it Exit times, moment problems
  and comparison theorems.} Potential Analysis {\bf 31} 1-8. 

\bibitem[MM1]{MM1} P. McDonald and R. Meyers (2003) {\it Dirichlet
  spectrum and heat content.} J. Funct. Anal. {\bf 200} 150--159.

\bibitem[MM2]{MM2} P. McDonald and R. Meyers (2003) {\it Isospectral
  polygons, planar graphs and heat content.} Proc. AMS {\bf 131} 3589--3599.

\bibitem[KPS1]{KPS1} V. Kostrykin, J. Potthoff and R. Schrader  (2007)
  {\it Heat kernels on metric graphs and a trace formula.}  Adventures
  in mathematical physics,  175--198, Contemp. Math., 447,
  Amer. Math. Soc., Providence, RI.

\bibitem[KPS2]{KPS2}  V. Kostrykin, J. Potthoff and R. Schrader
  (2008) {\it Contraction semigroups on metric graphs.}  Analysis on
  graphs and its applications,  423--458,   Proc. Sympos. Pure Math.,
  77, Amer. Math. Soc., Providence, RI  

\bibitem[KPS3]{KPS3}  V. Kostrykin, J. Potthoff and R. Schrader
  (2005) {\it Laplacians on metric graphs: eigenvalues resolvents and
  semigroups}  Proceeding of the Conference on Quantum Graphs and
  Their Applications, Amer. Math. Soc., Providence, RI  

\bibitem[HMP1]{HMP1}  A. Hurtado, S. Markorvsen and V.Palmer (2009)
  {\it Torsional rigidity of submanifolds with controlled geometry.}
  Mathem. Ann. {\bf 344}, 511--542.

\bibitem[HMP2]{HMP2}  A. Hurtado, S. Markorvsen and V.Palmer (2010)
  {\it Comparison     of exit moment spectra for extrinsic metric
    balls.} arXiv   1009.1257v1 [math.DG].

\bibitem[BE]{BE}  J. Bolte and S. Endres (2009)
  {\it The trace formula for quantum graphs with general self adjoint
    boundary conditions.} Ann. Henri Poincare {\bf 10}, 189--223.

\bibitem[SS]{SS} T. Schapira and U.Smilansky (2006) {\it
  Quantum graphs which sound the same} in: Non-linear Dynamics and
  Fundamental Interactions, F Khanna and D. Matrasulov (eds.), Springer,
  Berlin 17--29. 

\bibitem[So]{So} M. Solomyak (2002) {\it On the eigenvalue estimates
  for the weighted Laplacian on metric graphs.} Nonlinear problems in
  mathematical physics and related topics, I, Int. Math. Ser. (N.Y.)
  vol. 1, Kluwer/Plenum, New York, 327--347.

\bibitem[Ch]{Ch} S. J. Chapman (1995) {\it Drums that sound the same.}
  Am. Math. Monthly, 124--138.

\bibitem[Po1]{Po1} G. Pólya, {\it Torsional rigidity, principal
  frequency, electrostatic capacity and symmetrization,} Quarterly of
  Applied Math., 6 (1948), pp. 267, 277.   

\end{thebibliography}
\end{document}